# Extended framework of Hamilton's principle applied to Duffing oscillation

Jinkyu Kim[1], Hyeonseok Lee[2], Jinwon Shin[*]

[1]Research Professor, School of Architecture and Architectural Engineering, Hanyang University, 55 Hanyangdaehak-ro, Sangnok-gu, 15588, Kyeonggi-do, Korea.

[2]Director, Korea Construction Safety Technology, 17 SuyYeong-ro, 1079-3 GwangAn-dong, 665beon-gil, SuYeong-gu, Busan, 48243, Korea.

[*]Assistant professor, Department of Architectural Engineering, Catholic Kwandong University, 24, Beomil-ro 579beon-gil Gangneung-si, 25601, Gangwon-do, Korea. Email: jshin@cku.ac.kr, Tel: +82-33-649-7694, Fax :+82-33-645-8120

## Abstract

The paper begins with a novel variational formulation of Duffing equation using the extended framework of Hamilton's principle (EHP). This formulation properly accounts for initial conditions, and it recovers all the governing differential equations as its Euler-Lagrange equation. Thus, it provides elegant structure for the development of versatile temporal finite element methods. Herein, the simplest temporal finite element method is presented by adopting linear temporal shape functions. Numerical examples are included to verify and investigate performance of non-iterative algorithm in the developed method.

# 1. Introduction

Nonlinear oscillators have been a great interest and the Duffing oscillator is a well-known example [1-15]. In particular, the forced-damped Duffing oscillator is a seminal system for the study of chaos in nonlinear dynamics. Analytical methods including the traditional perturbation method, the homotopy perturbation method, the variational iteration method, the parameter-expansion method, and the He's frequency-amplitude formulation summarized in [16] are not applicable and only a few methods such as the Laplace decomposition [17] and the homotopy perturbation transform [18] are recently introduced to analytically solve this forced-damped Duffing oscillator.

For variational approach to nonlinear oscillations, some attempts have been made [19-27]. Ultimately, they applied the Hamilton's principle [28, 29] to nonlinear oscillator systems and presented either analytical solutions or numerical methods. For example, He [22] provided the energy balanced method. This method combines the Hamilton's principle and the semi-inverse method [30] to get analytical solutions of free vibration of nonlinear oscillator systems including higher-order Duffing oscillator. D'Acunto [24, 25] and Zhang [26] further developed this approach and established the He's variational method. Also, there were some numerical methods of Duffing equation stemming from variational integrators ([31, 32]- variational integrators are numerical integrators obtained from discretizing the Hamilton's principle). For example, Leok and Shingel [32] introduced the $2^{nd}$-order variational integrator for the unforced-undamped Duffing oscillator by adopting Hermite interpolation polynomials and the Euler-Maclaurin quadrature formula. Every effort seemed effective to account for nonlinear oscillation in their interests. However, it should be noted that the Hamilton's principle has critical weakness called end-points constraint (thorough review on this subject is available in the chapter 5 of [33]), which means that initial conditions cannot be properly considered in the Hamilton's principle. Rather, the Hamilton's principle utilizes temporal boundary conditions: it assumes that positions of a system at start and end of the time-interval are known. In fact, for initial value problems, both position and velocity of a system are known at the beginning and it is a main subject to find out how the system evolves from its beginning to the end of the time-interval. Thus, serious theoretical inconsistency exists in the original framework of Hamilton's principle.

To resolve end-points constraint in the Hamilton's principle, the extended framework of Hamilton's principle (EHP) was recently established for linear elastic and viscoplastic systems [34-37] and this has been applied to heat diffusion [38] and thermoelasticity [39]. Key ideas of the EHP are identified as use of mixed Lagrangian formulation [40-45] and new definition of the functional action. The EHP is simple and has extensive applicability as the original Hamilton's principle has been served broadly throughout mathematical physics and engineering. The current work focuses on its application to the forced-damped Duffing equation. The remainder of the paper is organized as follows. Section 2 presents variational formulation of Duffing equation in the context of EHP, where all the governing differential equations are recovered from the corresponding Euler-Lagrange equations along with proper use of initial conditions. Then, in Section 3, with employing linear shape functions over time domain, this formulation is discretized to produce the simplest temporal finite element method. Some numerical examples are considered in Section 4 to validate this new computational method. Finally, in Section 5, the present work is summarized and conclusions are drawn.

## 2. Variational formulation of Duffing equation

In this section, a new variational formulation of the forced Duffing oscillator stemming from the extended framework of Hamilton's principle (EHP) is presented. In this approach, the Duffing oscillator is regarded as combination of linear and nonlinear contribution of elasticity as described in Fig. 1.

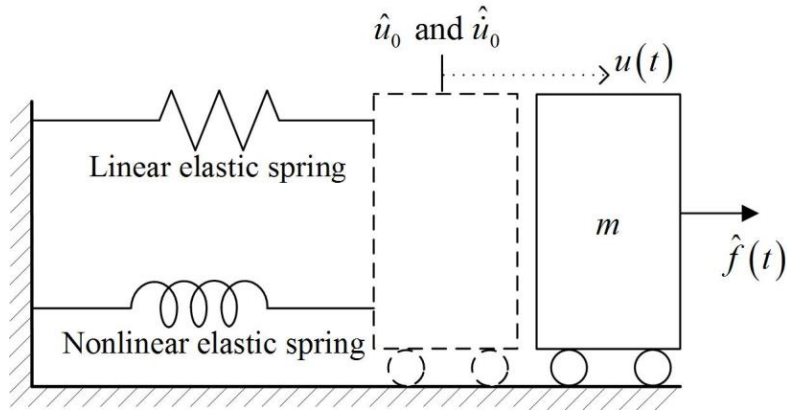

**Fig. 1. The Duffing oscillator in the EHP**

With mass $m$, linear elastic stiffness $k$, applied force $\hat{f}$, and effect of cubic nonlinearity, the forced-undamped Duffing equation is given by

$$m\ddot{u} + ku + \beta u^3 = \hat{f} \tag{1}$$

along with the initial conditions

$$u(0) = \hat{u}_0; \quad \dot{u}(0) = \hat{\dot{u}}_0 \tag{2}$$

In Eqs.(1)-(2), $u = u(t)$ represents a displacement at time $t$ and a superposed dot is a derivative with respect to time.

For this system, the EHP defines the Lagrangian $L$ in mixed variables as

$$L(u, \dot{u}, \dot{J}; t) = \frac{1}{2} m \, \dot{u}^2 + \frac{1}{2} a \, \dot{J}^2 - \dot{J} \, u - \frac{\beta}{4} u^4 + \hat{f} \, u \tag{3}$$

In Eq.(3), $a$ is the flexibility of linear elastic spring ($a = 1/k$) and $J(t)$ represents the impulse of linear elastic spring force $F_s(t)$. Thus,

$$J(t) = \int_0^t F_s(\tau) d\tau = \int_0^t k \, u(\tau) d\tau \tag{4}$$

For the time interval $[0, t]$, the functional action $I$ is written by

$$I = \int_0^t L(u, \dot{u}, \dot{J}; \tau) d\tau \tag{5}$$

Then, the EHP newly defines the first variation of the above action integral as

$$\delta I_{NEW} = -\delta \int_0^t L(u, \dot{u}, \dot{J}; \tau) d\tau + \underline{\left[ m \hat{\dot{u}} \, \delta \hat{u} \right]_0^t} \tag{6}$$

Eq.(6) consists of two parts. If just following the original Hamilton's principle, we only have the first term on the right-hand side in Eq.(6) $\left( \text{that is, } \delta I = -\delta \int_0^t L(u, \dot{u}, \dot{J}; \tau) d\tau \right)$.

However, in the EHP, the first variation of the action is newly defined by adding up the underlined terms. Such additional terms are the counterparts to the terms without end-points constraint in the original framework of Hamilton's principle.

Main objective of the present formulation is to show that the stationarity of the action $I$ adopting Eq.(6) rather than $\delta I$ can account for all the governing differential equations along with proper use of initial conditions of the forced-undamped Duffing oscillator depicted in Fig. 1.

After substituting Eq.(3) into Eq.(6), we have

$$\delta I_{NEW} = -\int_0^t \left[ m\,\dot{u}\,\delta\dot{u} + a\,\dot{J}\,\delta\dot{J} - u\,\delta\dot{J} - \dot{J}\,\delta u - \beta\,u^3\,\delta u + \hat{f}\,\delta u \right] d\tau + \left[ m\hat{\dot{u}}\,\delta\hat{u} \right]_0^t \quad (7)$$

Next, to remove all the time derivatives on variations, temporal integration by parts is performed on $m\dot{u}\,\delta\dot{u}$, $a\,\dot{J}\,\delta\dot{J}$, $-u\,\delta\dot{J}$ in Eq.(7). The resulting formulation becomes

$$\delta I_{NEW} = \underline{\left[ m\hat{\dot{u}}\,\delta\hat{u} \right]_0^t} - \left[ m\,\dot{u}\,\delta u \right]_0^t + \left[ \left( u - a\,\dot{J} \right)\,\delta J \right]_0^t$$
$$+ \int_0^t \left( m\,\ddot{u} + \dot{J} + \beta\,u^3 - \hat{f} \right)\delta u\,d\tau + \int_0^t \left( a\,\ddot{J} - \dot{u} \right)\delta J\,d\tau$$

$$(8)$$

Recalling the stationarity of the action as $\delta I_{NEW} = 0$ with arbitrary variations of $\delta u$ and $\delta J$ between the time interval leads to set of Euler-Lagrange equations as

$$m\,\ddot{u} + \dot{J} + \beta\,u^3 - \hat{f} = 0 \quad (9)$$

and

$$a\,\ddot{J} - \dot{u} = 0 \quad (10)$$

from the second line of Eq.(8), each of which represents dynamic equilibrium and rate-constitutive relation of linear elastic spring. Also, from the last term on the first line of Eq.(8), we have the linear elastic constitutive relation

$$u - a\,\dot{J} = 0 \quad (11)$$

at the ends of the time interval. Thus, in the EHP, all the governing differential equations are obtained in terms of displacement $u$ and impulsive linear elastic spring force $J$. More importantly, initial conditions are properly considered with the first two terms on the right-hand side fo Eq.(8) throughout the concept of sequentially assigning process. Notice that the unique dynamic evolution of the forced-undamped Duffing oscillator is guaranteed by enforcing unspecified unique values of $\dot{u}(0)$, $u(0)$, $\dot{u}(t)$, $\delta u(t)$ to $\hat{\dot{u}}(0)$, $\delta \hat{u}(0)$, $\hat{\dot{u}}(t)$, $\delta \hat{u}(t)$:

$$\hat{\dot{u}}(0) = \dot{u}(0); \quad \delta \hat{u}(0) = \delta u(0); \quad \hat{\dot{u}}(t) = \dot{u}(t); \quad \delta \hat{u}(t) = \delta u(t) \tag{12}$$

Now, let us first assign the specified initial velocity $\hat{\dot{u}}_0$ to $\hat{\dot{u}}(0)$ as

$$\hat{\dot{u}}(0) = \hat{\dot{u}}_0 \tag{13}$$

and sequentially assign the specified initial displacement $\hat{u}_0$ to $\hat{u}(0)$ as

$$\delta \hat{u}(0) = 0 \quad (\hat{u}(0) = \hat{u}_0) \tag{14}$$

The subsequent zero-valued term in Eq.(14) is not necessarily appeared on the newly defined action variation in Eq.(7). However, this supplementary conceptual condition is necessary for proper use of initial conditions in the EHP (the interested reader is referred to [34, 36]).

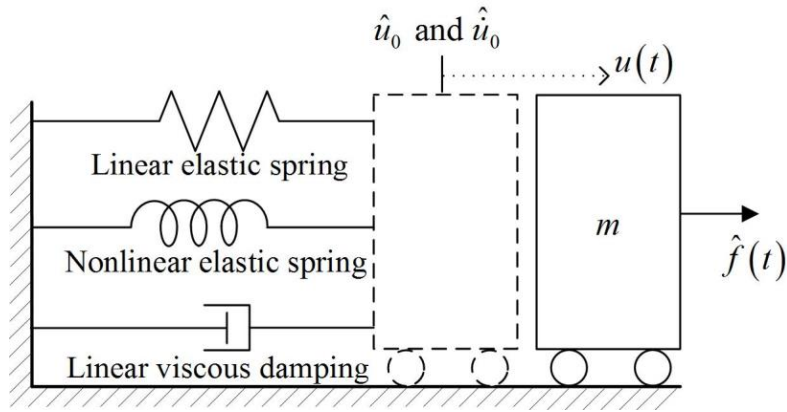

Fig. 2. The forced-damped Duffing oscillator

In the EHP, dissipation can be considered with the Rayleigh's formalism [46]. For the forced-damped Duffing oscillator as indicated in Fig. 2, the Rayleigh's dissipation $\varphi$ is given by

$$\varphi(\dot{u};\,t) = \tfrac{1}{2}c\big[\dot{u}(t)\big]^2 \tag{15}$$

and the subsequent variation of this Rayleigh's dissipation is defined as

$$\int_0^t \frac{\partial\varphi(\dot{u};\tau)}{\partial\dot{u}}\,\delta u\,d\tau \tag{16}$$

Thus, for this system, the action variation in the EHP is obtained by adding the following term to Eq.(7).

$$\int_0^t \frac{\partial\varphi(\dot{u};\tau)}{\partial\dot{u}}\,\delta u\,d\tau = \int_0^t c\,\dot{u}\,\,\delta u\,d\tau \tag{17}$$

Then, the resulting dynamic equilibrium equation becomes

$$m\,\ddot{u} + c\,\dot{u} + \dot{J} + \beta\,u^3 - \hat{f} = 0 \tag{18}$$

with leaving others to be the same as Eq.(10)-(14). Notice that with Eq.(4), the Eq.(18) can be equivalently written as

$$m\,\ddot{u} + c\,\dot{u} + k\,u + \beta\,u^3 - \hat{f} = 0 \tag{19}$$

that can be found in the displacement-based approach.

## 3. Temporal finite element method

The EHP is not a complete variational principle, since it requires the Rayleigh's formalism for dissipative process and the new action variation is not entirely satisfactory in variational sense. Even with such weakness, the EHP has simple structure and broad applicability as the original Hamilton's principle. Also, with proper use of initial conditions, it provides sound base for the development of novel computational methods involving finite element method over time.

In this section, we present probably the simplest temporal finite element method for the forced-damped Duffing oscillator. For this development, the following weak form is utilized.

$$\delta I_{NEW} = -\int_0^{\Delta t}\Big[ m\,\dot{u}\,\delta\ddot{u} + c\,\dot{u}\,\delta u + a\,\dot{J}\,\delta\dot{J} - u\,\delta\dot{J} - \dot{J}\,\delta u - \beta u^3\,\delta u + \hat{f}\,\delta u \Big]d\tau$$
$$+\Big[\hat{p}\,\delta\hat{u}\Big]_0^{\Delta t} = 0 \tag{20}$$

Here, we consider a small time-duration from 0 to $\Delta t$ and separately introduce the linear momentum $\hat{p}\left(=m\,\hat{\dot{u}}\right)$ coming from the additional terms in Eq.(6). Main reason to have $\hat{p}$ at the ends of the time-duration is for the explicit use of initial velocity in a numerical method.

Eq. (20) has continuity-balanced variables. That is, primary variables $u(t)$ and $J(t)$ as well as the variations $\delta u(t)$ and $\delta J(t)$ have $C^0$ temporal continuity requirement, while the linear momentum $\hat{p}$ has $C^{-1}$ temporal continuity requirement. Thus, $\hat{p}$ can be approximated as constant at the ends of the time-duration such as $\hat{p}_0$ and $\hat{p}_{\Delta t}$. Also, one can approximate $u$, $\delta u$, $J$, and $\delta J$ with following linear shape functions on time domain.

$$L_0(\tau) = 1 - \frac{\tau}{\Delta t} \tag{21}$$

$$L_1(\tau) = \frac{\tau}{\Delta t} \tag{22}$$

Subsequently, $\dot{u}$, $\delta\dot{u}$, $\dot{J}$, and $\delta\dot{J}$ can be approximated with first derivatives of above shape functions as

$$\dot{L}_0(\tau) = -\frac{1}{\Delta t} \qquad (23)$$

$$\dot{L}_1(\tau) = \frac{1}{\Delta t} \qquad (24)$$

With Eqs.(21)-(24), for a representative case, let us discretize the first term in Eq.(20) as follows.

$$
\begin{aligned}
-\int_0^{\Delta t} m\, \dot{u}\, \delta\dot{u}\, d\tau &= -\int_0^{\Delta t} m\, \delta\dot{u}\, \dot{u}\, d\tau \\
&= -m \int_0^{\Delta t} \lfloor \delta u_0 \quad \delta u_{\Delta t} \rfloor \left\{ \begin{matrix} \dot{L}_0 \\ \dot{L}_1 \end{matrix} \right\} \lfloor \dot{L}_0 \quad \dot{L}_1 \rfloor \left\{ \begin{matrix} u_0 \\ u_{\Delta t} \end{matrix} \right\} d\tau \\
&= -m \lfloor \delta u_0 \quad \delta u_{\Delta t} \rfloor \int_0^{\Delta t} \begin{bmatrix} \dot{L}_0\, \dot{L}_0 & \dot{L}_0\, \dot{L}_1 \\ \dot{L}_1\, \dot{L}_0 & \dot{L}_1\, \dot{L}_1 \end{bmatrix} d\tau \left\{ \begin{matrix} u_0 \\ u_{\Delta t} \end{matrix} \right\} \\
&= -m \lfloor \delta u_0 \quad \delta u_{\Delta t} \rfloor \begin{bmatrix} \dfrac{1}{\Delta t} & -\dfrac{1}{\Delta t} \\ -\dfrac{1}{\Delta t} & \dfrac{1}{\Delta t} \end{bmatrix} \left\{ \begin{matrix} u_0 \\ u_{\Delta t} \end{matrix} \right\} \\
&= \left( -\frac{m}{\Delta t} u_0 + \frac{m}{\Delta t} u_{\Delta t} \right) \delta u_0 + \left( \frac{m}{\Delta t} u_0 - \frac{m}{\Delta t} u_{\Delta t} \right) \delta u_{\Delta t}
\end{aligned}
\qquad (25)
$$

Here, first, we introduce discretized representations for $\dot{u}$ and $\delta\dot{u}$ with first derivatives of shape functions and discrete variables $u_0$, $u_{\Delta t}$, $\delta u_0$, $\delta u_{\Delta t}$. Then, integral computation is undertaken. Lastly, the resulting formulation is expressed separately with respect to variations such as $\delta u_0$ and $\delta u_{\Delta t}$.

After applying this step-by-step approach to all the other terms and collecting terms with respect to relevant variations, we have

$$\delta I_{NEW} = \left[ \frac{m}{\Delta t}(u_{\Delta t} - u_0) + \frac{c}{2}(u_{\Delta t} - u_0) + \frac{1}{2}(J_{\Delta t} - J_0) - \hat{f}_{\Delta t}\frac{\Delta t}{2} - \hat{p}_0 + Q \right]\delta u_0$$

$$+ \left[ -\frac{m}{\Delta t}(u_{\Delta t} - u_0) + \frac{c}{2}(u_{\Delta t} - u_0) + \frac{1}{2}(J_{\Delta t} - J_0) - \hat{f}_{\Delta t}\frac{\Delta t}{2} + \hat{p}_{\Delta t} + Q \right]\delta u_{\Delta t}$$

$$+ \left[ \frac{a}{\Delta t}(J_{\Delta t} - J_0) - \frac{1}{2}(u_{\Delta t} + u_0) \right]\delta J_0 \tag{26}$$

$$+ \left[ -\frac{a}{\Delta t}(J_{\Delta t} - J_0) + \frac{1}{2}(u_{\Delta t} + u_0) \right]\delta J_{\Delta t} = 0$$

with $Q$ representing a cubic function of $u_0$ and $u_{\Delta t}$ given by

$$Q(u_0, u_{\Delta t}) = \frac{\beta \Delta t}{20}\left( 4u_0^3 + 3u_{\Delta t}u_0^2 + 2u_{\Delta t}^2 u_0 + u_{\Delta t}^3 \right) \tag{27}$$

While deriving Eq.(26), the external forcing $\hat{f}$ is assumed to be constant over the time-duration $[0, t]$ with values of $\hat{f}_{\Delta t}$. Also, the relation $\delta u = \delta\hat{u}$ at the ends of the time-duration if used with the enforced condition in Eq.(12).

Previously, in the EHP, the concept of sequentially assigning initial values is utilized, where subsequently imposed initial displacement is not explicitly appeared. Following this idea strictly, to satisfy the discretized representation of stationarity of the action in Eq.(26), each coefficient of discrete variations must be taken as zero. Consequently, one can obtain the following independent equations from Eq.(26)- notice that the third row and the fourth row in Eq.(26) result in one equation.

$$\frac{m}{\Delta t}(u_{\Delta t} - u_0) + \frac{c}{2}(u_{\Delta t} - u_0) + \frac{1}{2}(J_{\Delta t} - J_0) - \hat{f}_{\Delta t}\frac{\Delta t}{2} - \hat{p}_0 + Q = 0 \tag{28}$$

$$-\frac{m}{\Delta t}(u_{\Delta t} - u_0) + \frac{c}{2}(u_{\Delta t} - u_0) + \frac{1}{2}(J_{\Delta t} - J_0) - \hat{f}_{\Delta t}\frac{\Delta t}{2} + \hat{p}_{\Delta t} + Q = 0 \tag{29}$$

$$\frac{a}{\Delta t}(J_{\Delta t} - J_0) - \frac{1}{2}(u_{\Delta t} + u_0) = 0 \tag{30}$$

where

$$Q = \frac{\beta \Delta t}{20} \left( 4 u_0^3 + 3 u_0^2 u_{\Delta t} + 2 u_0 u_{\Delta t}^2 + u_{\Delta t}^3 \right) \tag{31}$$

With use of MAPLE13 [47], above equations can be solved as

$$u_{\Delta t} = RootOf\left( \bar{Q} \right) \tag{32}$$

$$J_{\Delta t} = J_0 + \frac{\Delta t}{2 a} \left( u_0 + u_{\Delta t} \right) \tag{33}$$

$$\hat{p}_{\Delta t} = \frac{1}{4 \Delta t\, a} \begin{bmatrix} \beta\, h^2\, a\, u_0 (u_{\Delta t}^2 + 2 u_0 u_{\Delta t} + 3 u_0^2) + 3 \Delta t^2 \left( u_{\Delta t} + u_0 \right) \\ + 20\, a\, m \left( u_{\Delta t} - u_0 \right) + 6\, c\, a\, \Delta t \left( u_{\Delta t} - u_0 \right) \\ + 16\, a\, \Delta t\, \hat{p}_0 - 6\, a\, \Delta t^2\, \hat{f}_{\Delta t} \end{bmatrix} \tag{34}$$

Eq.(32) represents that $u_{\Delta t}$ is a solution of the following cubic equation $\bar{Q}$

$$\bar{Q} = A x^3 + B x^2 + C x + D = 0 \tag{35}$$

with coefficients

$$A = \beta\, a\, \Delta t^2 \tag{36}$$

$$B = 2 \beta\, a\, \Delta t^2\, u_0 \tag{37}$$

$$C = 5 \Delta t^2 + 20\, m\, a + 10\, c\, a\, \Delta t + 3 \beta\, a\, u_0^2\, \Delta t^2 \tag{38}$$

$$D = 5 u_0 \Delta t^2 - 10\, c\, a\, u_0\, \Delta t - 20\, a\, \hat{p}_0\, \Delta t - 20\, m\, a\, u_0 + 4 \beta\, a\, u_0^3\, \Delta t^2 - 10\, a\, \hat{f}_{\Delta t}\, \Delta t^2 \tag{39}$$

Considering physical interpretation of $u_{\Delta t}$ (the displacement at $\Delta t$) as well as characteristics of a cubic equation (existence of at least one real-valued solution), we expect $u_{\Delta t}$ to have a real number. This only happens when Eq.(35) has one real-valued root- if $u_{\Delta t}$ is a real number solution, then, correspondingly, one can expect $J_{\Delta t}$ and $\hat{p}_{\Delta t}$ in Eqs.(33)-(34) to have a real number solution, respectively. However, it would be problematic if Eq.(35) has more than two real-valued roots. Criterion for a cubic equation to have one real-valued root is given by the following statement.

For $x^3 + A_1 x^2 + B_1 x + C_1 = 0$, $\Delta$ must be greater than $0$ $(\Delta > 0)$. (40)

Here $\Delta$ is the discriminant of a cubic equation given by

$$\Delta = \frac{r^2}{4} + \frac{s^3}{27} \tag{41}$$

with

$$r = \frac{2A_1^3}{27} - \frac{A_1 B_1}{3} + C_1; \quad s = B_1 - \frac{A_1^3}{3} \tag{42}$$

In Eqs.(40)-(42), one can notice that if $s$ is greater than $0$, then, the cubic equation has always only one real-valued root, which is desirable to solve Eq.(35). After dividing Eq.(35) with respect to $A$ in Eq.(36), identifying $s$ from Eq.(42), and invoking $s > 0$, we have

$$\left( \frac{15\Delta t^2 + 10 c\,a\,\Delta t + 20 m\,a + 3 a\,u_0^2}{\beta\,a\,\Delta t^2} \right) - 8u_0^3 > 0 \tag{43}$$

Eq.(43) is the condition for $\Delta t$ to have unique real-valued solution of $u_{\Delta t}, J_{\Delta t}, \hat{p}_{\Delta t}$ in the developed numerical method. Thus, one can think that the time-step $\Delta t$ satisfying Eq.(43) adjusts Eq.(35) to have one-real valued root. The only concern is about nonexistence of a real-valued positive $\Delta t$ satisfying Eq.(43) and the decision making process on $\Delta t$. In numerical analysis, usually, with consideration of accuracy, a small positive $\Delta t$ is fixed before any computation begins. If such strategy is followed in the present method, the fixed real-valued $\Delta t$ may not satisfy Eq.(43) at every time-step, and in turn, Eq.(35) may yield more than one real-valued roots. The safety gadget to fix this situation requires the refined consideration on quality of $u_{\Delta t}$:

Among real-valued solutions of $u_{\Delta t}$, the true $u_{\Delta t}$ makes $\left| u_{\Delta t} - u_0 \right|$ smallest.     (44)

The statement (44) is quite reasonable since one can expect small changes in $u_{\Delta t}$ when taking small $\Delta t$ in numerical analysis.

In summary, the present temporal finite element method proceeds as follows.

Step 1. Identify system parameters and initial conditions ($m$ ; $a$ ; $\beta$ ; $c$ ; $\hat{f}_{\Delta t}$ ; $u_0$ ; $J_0$ ; $\hat{p}_0$).

Step 2. Identify temporal parameters (time-step $\Delta t$ and total analysis time $t$ )

Step 3. Solve $u_{\Delta t}$ with Eq.(32), Eqs.(36)-(39), and the statement (44).

Step 4. Solve $J_{\Delta t}$ and $\hat{p}_{\Delta t}$ with Eqs.(33)-(34).

Step 5. Store $u_{\Delta t}$ , $J_{\Delta t}$ , $\hat{p}_{\Delta t}$ .

Step 6. Update $u_{\Delta t}$ , $J_{\Delta t}$ , $\hat{p}_{\Delta t}$ (these become $u_0$ , $J_0$ , $\hat{p}_0$ in the next step).

Step 7. Repeat Step3~Step 6 until the total analysis time $t$ .

We would like to clarify two points in the above step-by-step algorithm.

First, the initial value $J_0$ in the step 1 is taken as zero. Let us consider the following equation.

$$m\dot{u}(t) + cu(t) + J(t) + \int_{t_p}^{t} \beta\left[u(\tau)\right]^3 d\tau - \int_{t_p}^{t} \hat{f}(\tau)d\tau + \bar{C}_1 = 0 \qquad (45)$$

In Eq.(45), $t_p$ represents a certain temporal point ($t_p < t$) and $\bar{C}_1$ is an arbitrary constant. Let us differentiate Eq.(45) with respect to $t$ by following the Leibniz integral rule. The resulting equation is exactly the same as Eq.(9). Thus, Eq.(45) holds at any time point $t$, where $\bar{C}_1$ remains unspecified. In fact, Eq.(45) can also be obtained from integrating Eq.(18) over $[t_p, t]$. This time, $\bar{C}_1$ is written as $\bar{C}_1 = -m\dot{u}(t_p) - cu(t_p) - J(t_p)$, where $\dot{u}(t_p)$, $u(t_p)$, and $J(t_p)$ are unspecified. In either way, $C_1$ cannot be specified.

Letting $t = 0$ in Eq.(45) results in

$$m\dot{u}_0 + cu_0 + J_0 + \int_{t_p}^{0} \beta[u(\tau)]^3 d\tau - \int_{t_p}^{0} \hat{f}(\tau) d\tau + \bar{C}_1 = 0 \qquad (46)$$

On the left-hand side of Eq.(46), the first two terms are evaluated with given initial conditions, while the last three terms cannot be decided. Therefore, $J_0$ can be any value. In fact, physically meaningful variable is $\dot{J}$ (not $J$) as shown in Eqs.(8)-(11) and Eq.(18), where the linear elastic spring force $\dot{J}$ and its derivative are involved. This can also be checked in the discrete version of stationarity of the action in the EHP: in Eqs.(28)-(31), $\dot{J}$ is approximated as $(J_{\Delta t} - J_0)/\Delta t$ with $J_0$ providing a reference. This is why we simply take $J_0 = 0$.

We would like to conclude this section with pointing out that unlike other numerical methods [48, 49], the present method is non-iterative to analyze the forced-damped Duffing equation.

## 4. Computational examples

The present numerical method is verified for three types of the forced-damped Duffing equation (19): Type I called "hard spring Duffing oscillator" ($k > 0$ and $\beta > 0$); Type II called "soft spring Duffing oscillator" ($k > 0$ and $\beta < 0$); Type III called "inverted Duffing oscillator" ($k < 0$ and $\beta > 0$). The "nonharmonic Duffing oscillator ( $k = 0$ and $\beta > 0$ in Eq.(19))" or the "Ueda oscillator" is not tested here because theoretically this cannot be accommodated within the current approach due to the existence of linear elastic spring force [see Fig. 1 and Eq.(18)].

For each type, we provide relevant examples found in [50-52]. Numerical simulation results obtained from the present non-iterative method are compared to those obtained from the Runge-Kutta-Fehlberg method (RKF45) using MAPLE13 [47]. Default values associated with iterations and convergence criteria are adopted for RKF45. For example, the default values for abserr (absolute error tolerance) and relerr (relative error tolerance) are given by 1e-7 and 1e-6, respectively.

In simulation tests, the time-step $\Delta t = 0.01$ is taken as a default value for the present method. However, we find that $\Delta t = 0.01$ is not sufficiently small for two examples in Type III (T3E5 and T3E6) to catch proper responses. Thus, for those examples, $\Delta t = 0.001$ is adopted. After method-comparison study, further, we investigate convergence characteristics of the present method with additional simulation for the representative example in each type.

## 4.1. Type I (hard spring Duffing oscillator)

For the Type I simulation, we take four examples found in [52], where system parameters in Eq.(18) are fixed as $m = 1$, $c = 0.2$, $k = 1$, $\beta = 0.1$, and $\hat{f} = 0.5\cos(2.00649t)$. Thus, here, initial conditions are control parameters and these are summarized in Table 1.

**Table 1. Control parameters for the Type 1 simulation**

| No. | Control parameters |
| --- | --- |
| T1E1 | $u_0 = 3$ and $\dot{u}_0 = 0$ |
| T1E2 | $u_0 = -3$ and $\dot{u}_0 = 0$ |
| T1E3 | $u_0 = -1$ and $\dot{u}_0 = 1$ |
| T1E4 | $u_0 = 1$ and $\dot{u}_0 = 1$ |

Fig. 3-Fig. 6 show numerical results of displacement-history and phase portraits for different initial conditions. As shown in all the figures, the present method yields almost identical results to RKF45 for the hard spring Duffing oscillator.

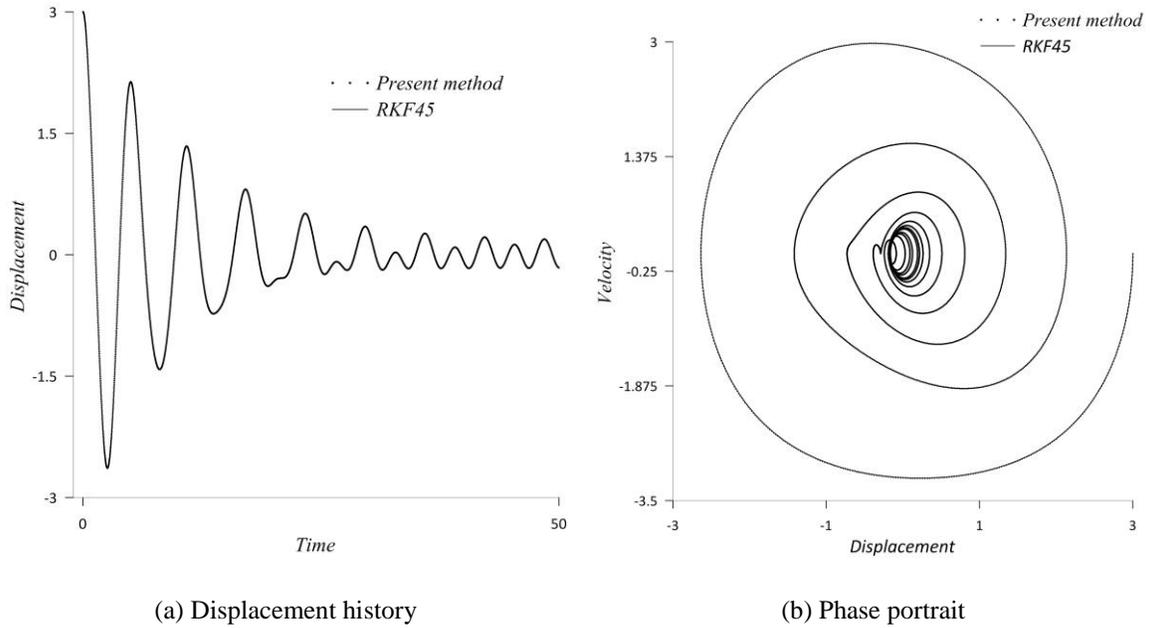

(a) Displacement history        (b) Phase portrait

**Fig. 3. Simulation results for T1E1**

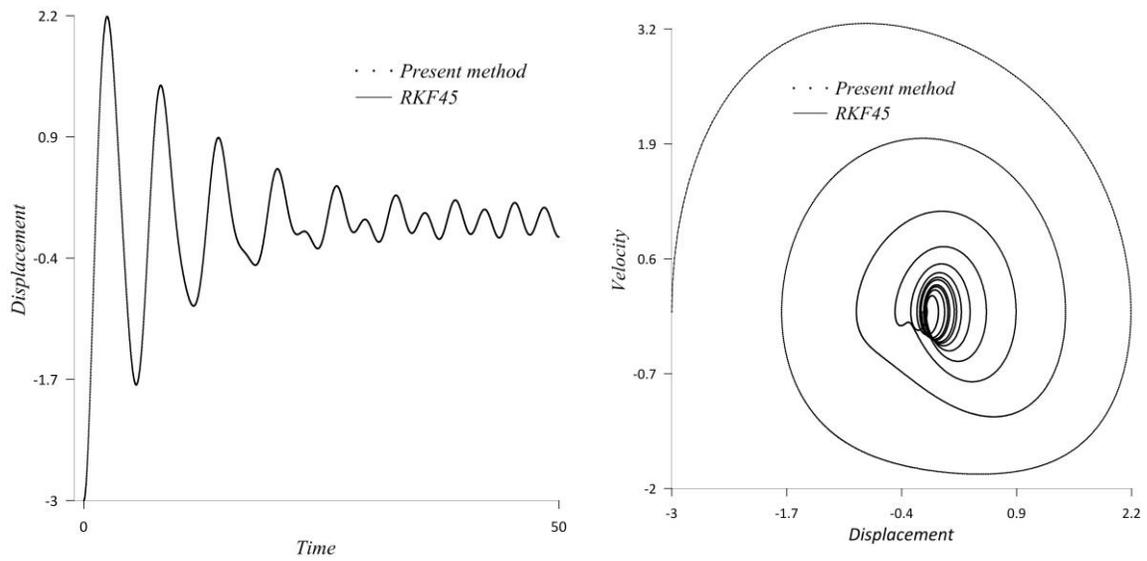

(a) Displacement history           (b) Phase portrait

**Fig. 4. Simulation results for T1E2**

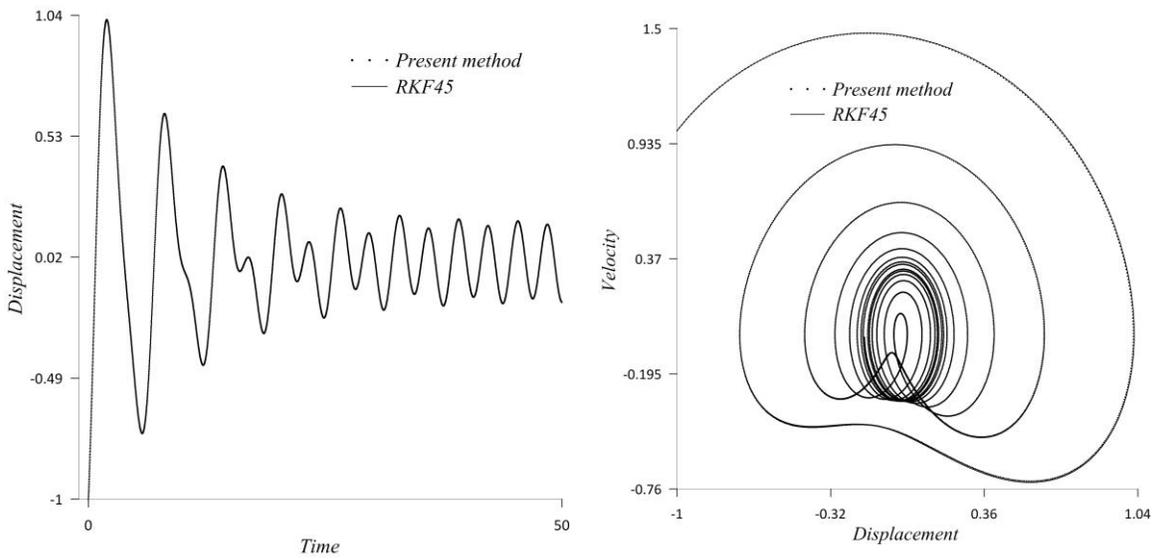

(a) Displacement history           (b) Phase portrait

**Fig. 5. Simulation results for T1E3**

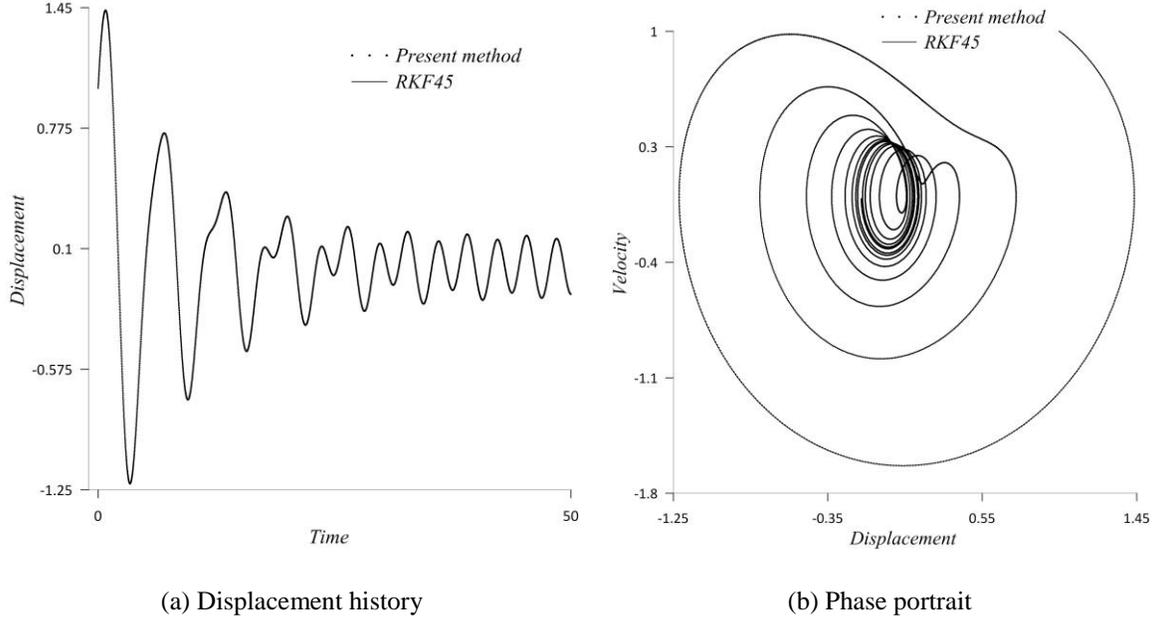

(a) Displacement history            (b) Phase portrait

**Fig. 6. Simulation results for T1E4**

## 4.2. Type II (soft spring Duffing oscillator)

Four examples excerpted from [51] are provided for the Type II simulation, where fixed system parameters in Eq.(19) are $m = 1$, $k = 1$, $\beta = -1/6$, and $\hat{f} = 1/3\cos(0.6t)$, while control parameters such as $c$ and initial conditions are specified in Table 2.

**Table 2. Control parameters for the Type II simulation**

| No. | Control parameters |
|---|---|
| T2E1 | $c = 0.24$, $u_0 = 0.519674$, $\dot{u}_0 = 0.072267$ |
| T2E2 | $c = 0.24$, $u_0 = 1$, $\dot{u}_0 = 0$ |
| T2E3 | $c = 0.002$, $u_0 = 0.55404958$, $\dot{u}_0 = 0.0011051$ |
| T2E4 | $c = 0.002$, $u_0 = 1$, $\dot{u}_0 = -0.531$ |

All the results are presented in Fig. 7-Fig. 10, where differences between RKF45 and the present method can hardly be detected.

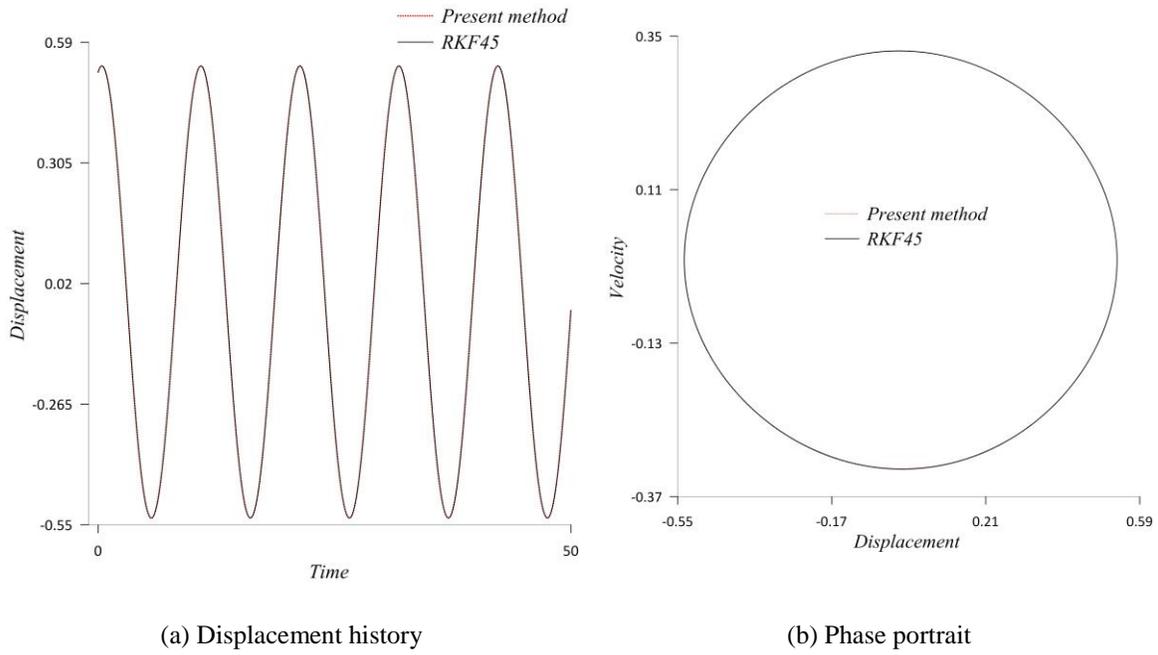

(a) Displacement history　　　　　　　　　(b) Phase portrait

**Fig. 7. Simulation results for T2E1**

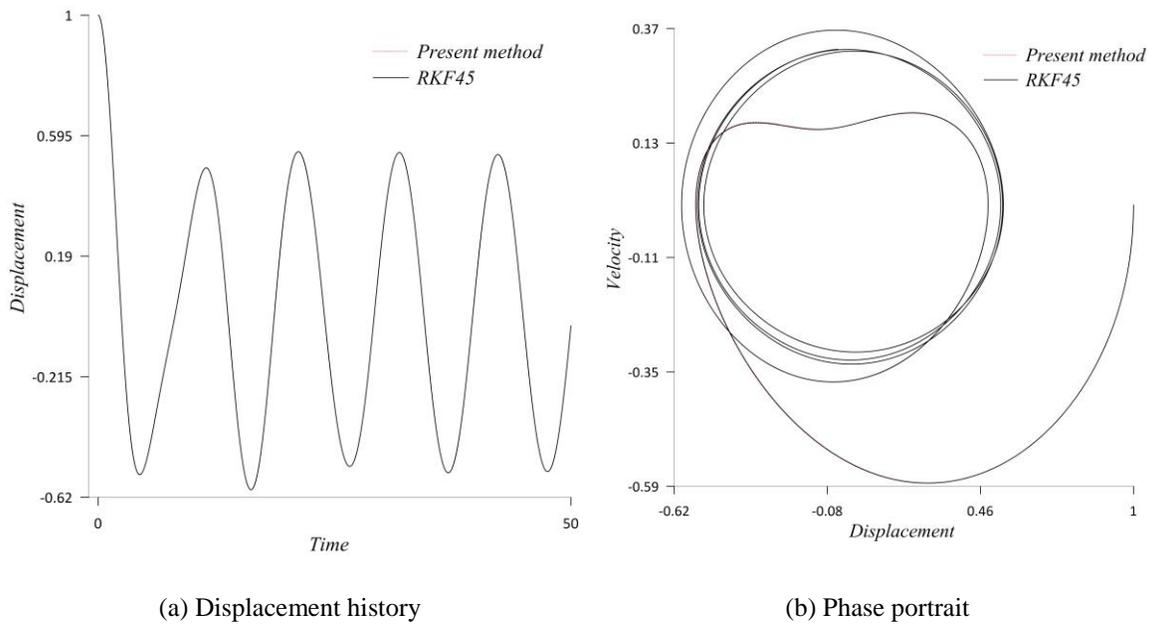

(a) Displacement history　　　　　　　　　(b) Phase portrait

**Fig. 8. Simulation results for T2E2**

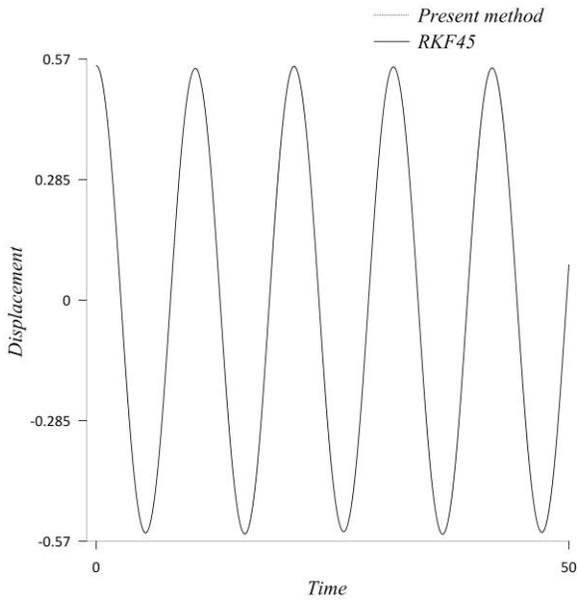

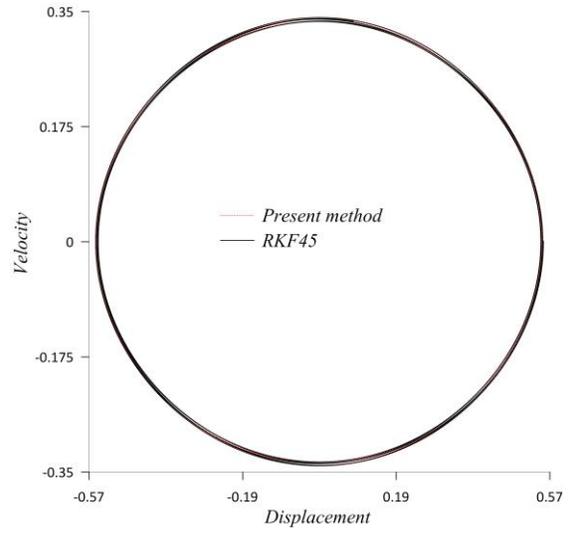

(a) Displacement history

(b) Phase portrait

**Fig. 9. Simulation results for T2E3**

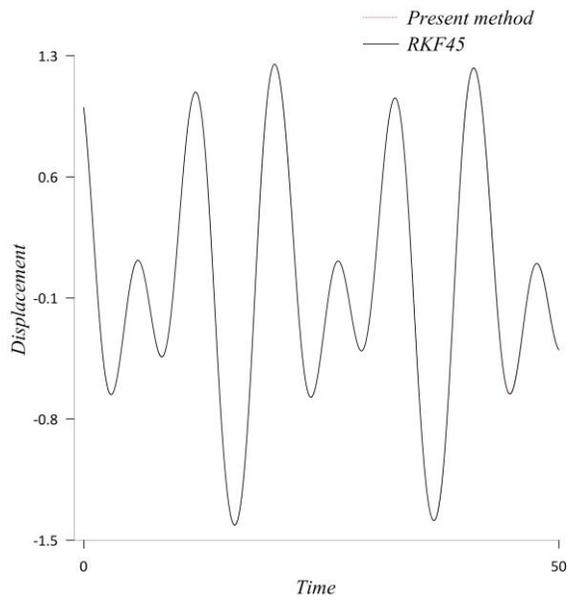

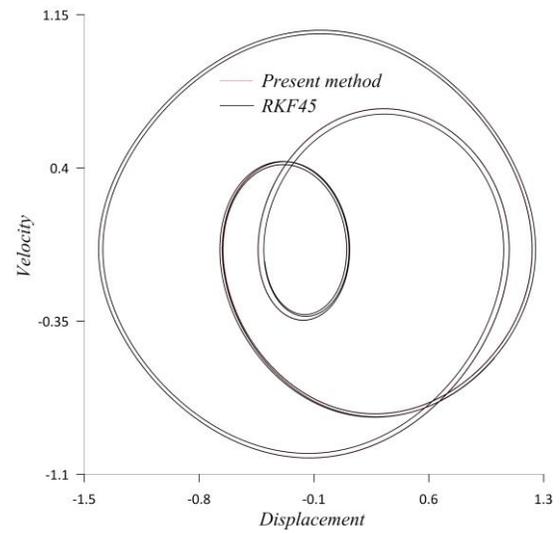

(a) Displacement history

(b) Phase portrait

**Fig. 10. Simulation results for T2E4**

*4.3. Type III (inverted Duffing oscillator)*

A total of six examples are provided for the Type III simulation. These examples were previously taken in [50] to show a great variety of periodic and chaotic solutions in the forced Duffing equation. System parameters and initial conditions are fixed as $m = 1$, $c = 0.3$, $k = -1$, $\beta = 1$, $u_0 = 1$, and $\dot{u}_0 = 0$. Control parameter is the amplitude of applied force $F_0$ given in $\hat{f} = F_0 \cos(1.2t)$, and this is specified in Table 3.

**Table 3. Control parameter for the Type III simulation**

| No. | Control parameter |
|---|---|
| T3E1 | $F_0 = 0.2$ |
| T3E2 | $F_0 = 0.28$ |
| T3E3 | $F_0 = 0.29$ |
| T3E4 | $F_0 = 0.37$ |
| T3E5 | $F_0 = 0.5$ |
| T3E6 | $F_0 = 0.65$ |

As shown in Fig. 11-Fig. 16, depending on different values of $F_0$, one can check the appearance of subharmonics through period-doubling bifurcation as well chaotic behavior. Also, in every figure, one can check that the present method shows good agreement to RKF45 for the inverted Duffing oscillator.

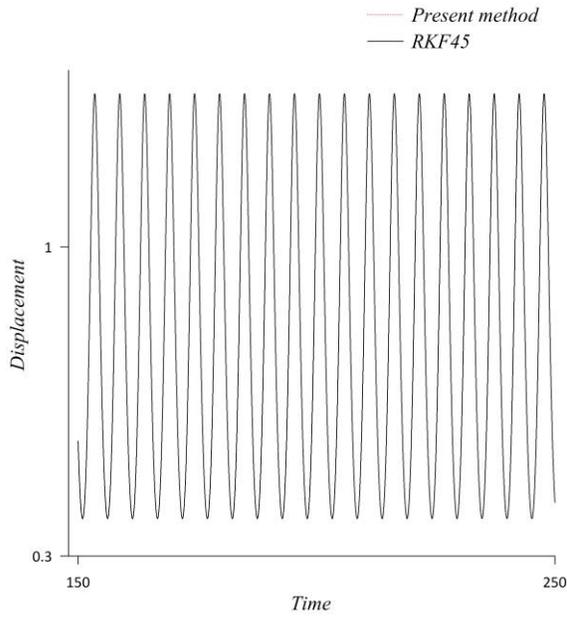
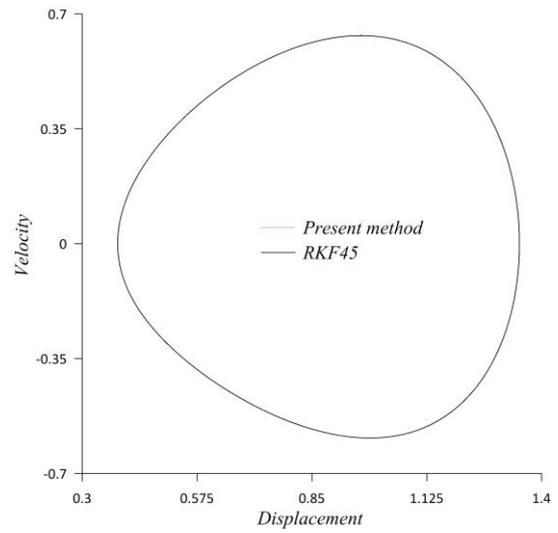

(a) Displacement history

(b) Phase portrait

**Fig. 11. Simulation results for T3E1**

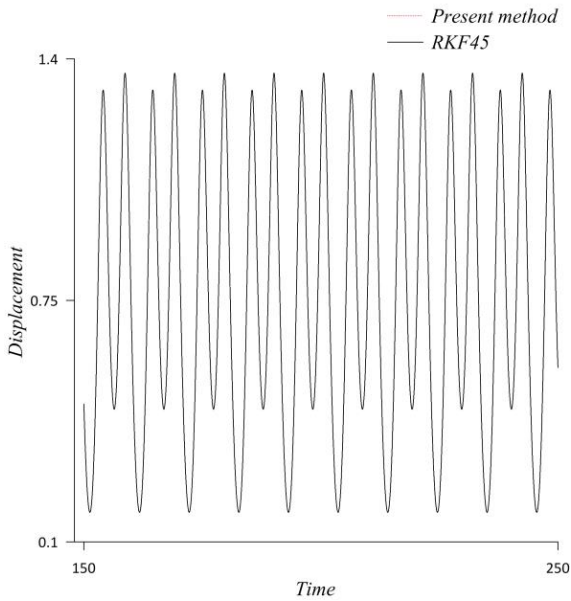
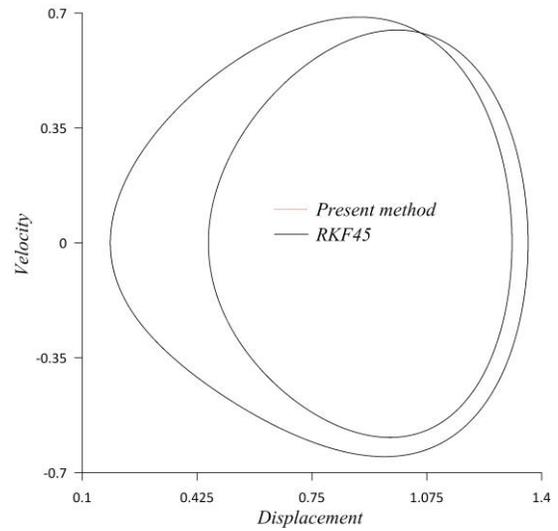

(a) Displacement history

(b) Phase portrait

**Fig. 12. Simulation results for T3E2**

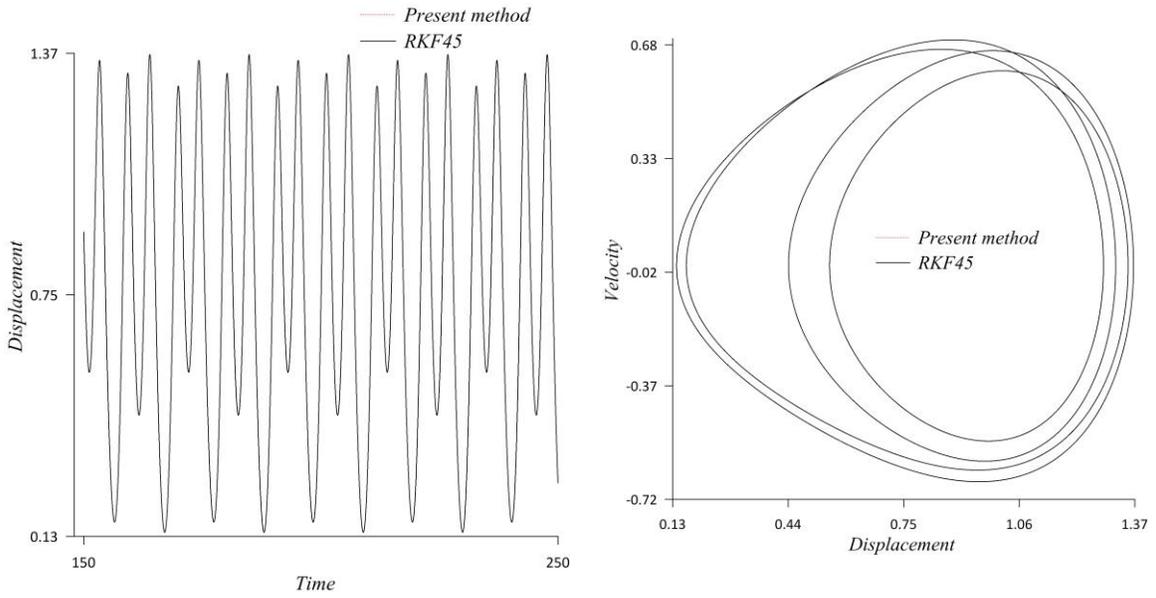

(a) Displacement history            (b) Phase portrait

**Fig. 13. Simulation results for T3E3**

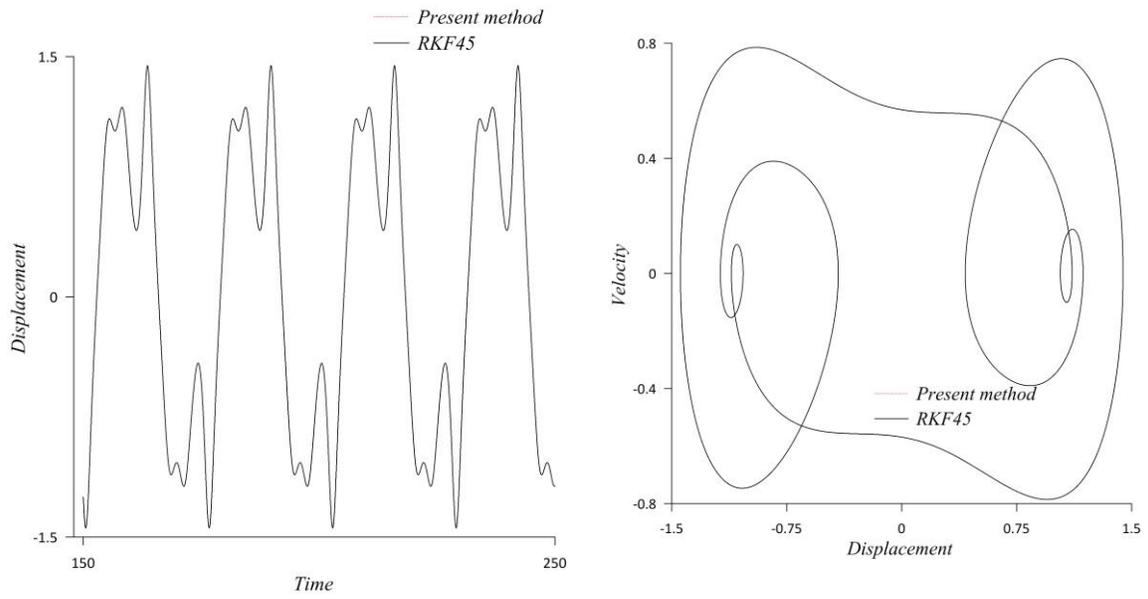

(a) Displacement history            (b) Phase portrait

**Fig. 14. Simulation results for T3E4**

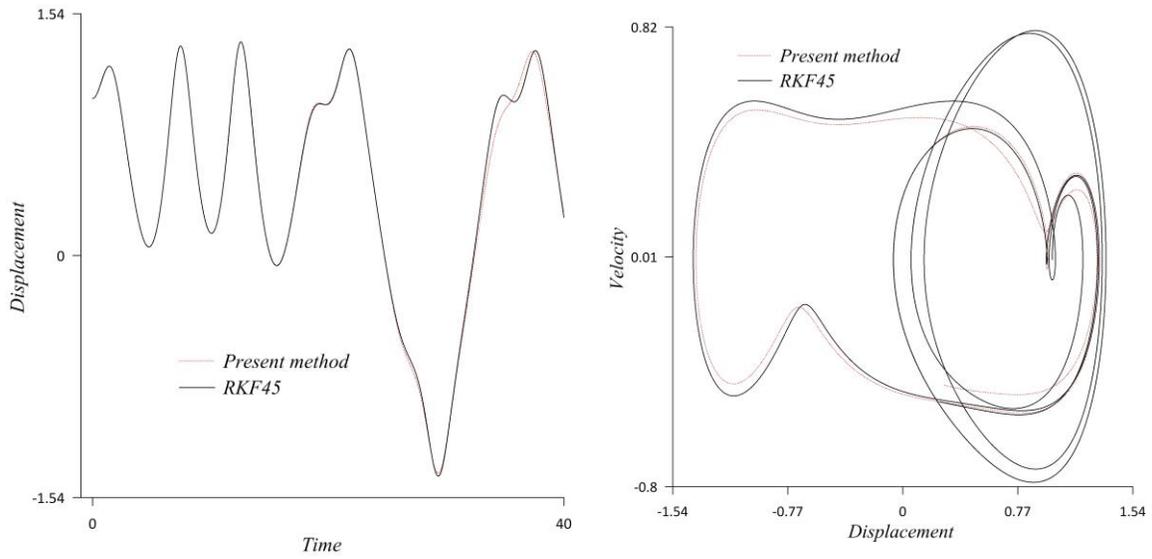

(a) Displacement history        (b) Phase portrait

**Fig. 15. Simulation results for T3E5**

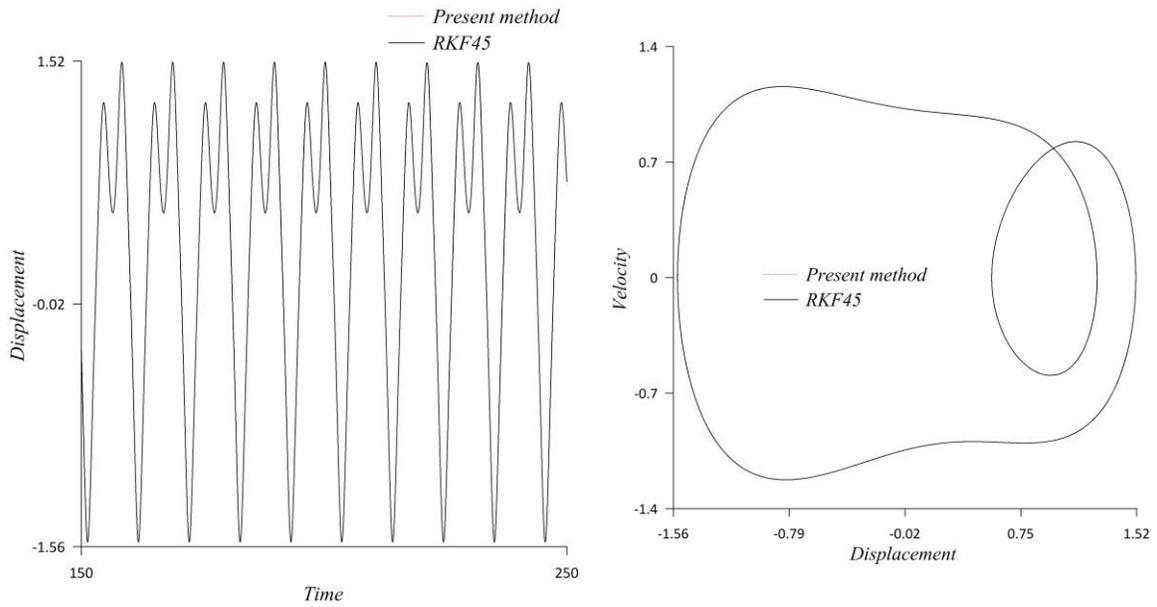

(a) Displacement history        (b) Phase portrait

**Fig. 16. Simulation results for T3E6**

Compared to other examples, relatively small analysis time $t \in [0, 40]$ is taken for

analyzing T3E5, and results between the present method and RKF45 show some

differences. Main reason for taking $t \in [0, 40]$ is that sufficient accuracy for RKF45 in

MAPLE with default values of iterations and convergence criteria resides there. Fig. 17

compares simulation results of RKF45 in MAPLE for the T3E5 problem, when different values of abserr and relerr are adopted.

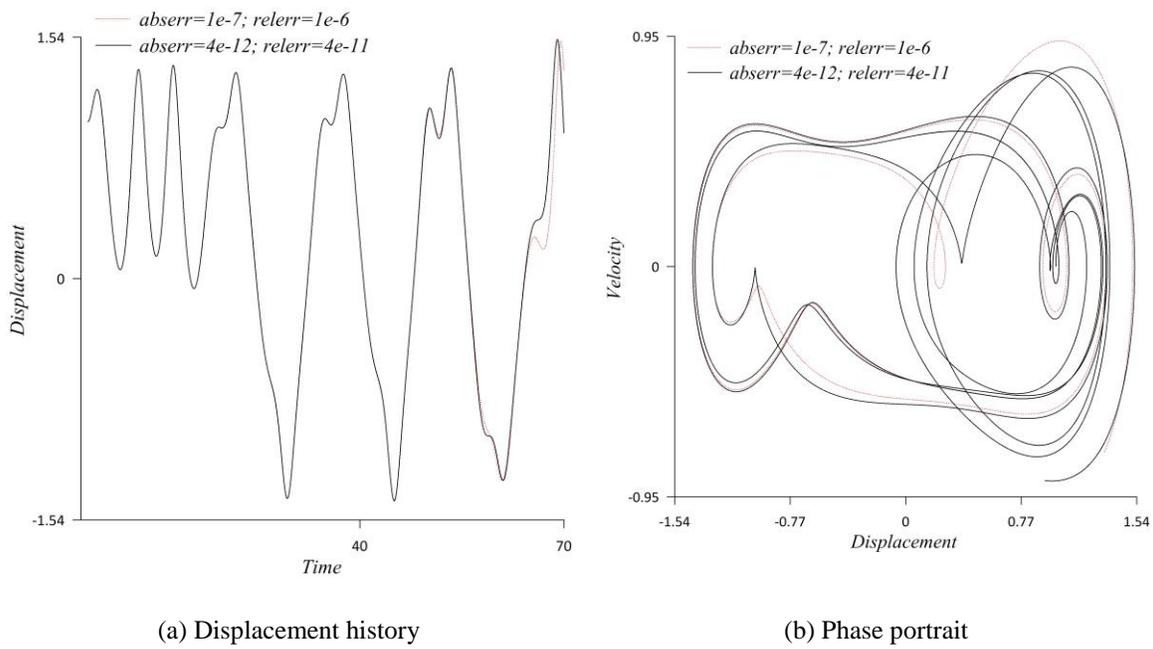

(a) Displacement history                     (b) Phase portrait

**Fig. 17. Simulation results for T3E5 obtained from RKF45 using different error criteria**

As shown in Fig. 17, one can see that accumulation of errors deteriorates the last of

analysis when using default values of abserr (1e-7) and relerr (1e-6) in RKF45, especially

for phase portrait. Thus, without specifying values of iteration and error criteria in

RKF45, one cannot reproduce the graph of T3E5 given in [50].

*4.4. Convergence characteristics in numerical methods*

Previously, for three types of the Duffing oscillation with relevant examples, we showed that the present non-iterative method using small time-steps ($\Delta t = 0.01$ for T1E1~T3E4 and $\Delta t = 0.001$ for T3E5~T3E6) and RKF45 in MAPLE using default values for iterations and convergence criteria are comparable to each other. Here, we investigate convergence characteristics of the present method. For this purpose, we take T1E1, T2E4, and T3E1 as the representative case for each type of Duffing oscillation. Three different values of $\Delta t$ ($\Delta t = 0.5$; $\Delta t = 0.2$; $\Delta t = 0.1$) are considered with fixing the analysis time as $t = 100$.

Fig. 18-Fig. 20 show a full trajectory of representative case on the phase portrait, where results are presented in two sub-regions such as $t \in [0, 30]$ and $t \in [30, 100]$.

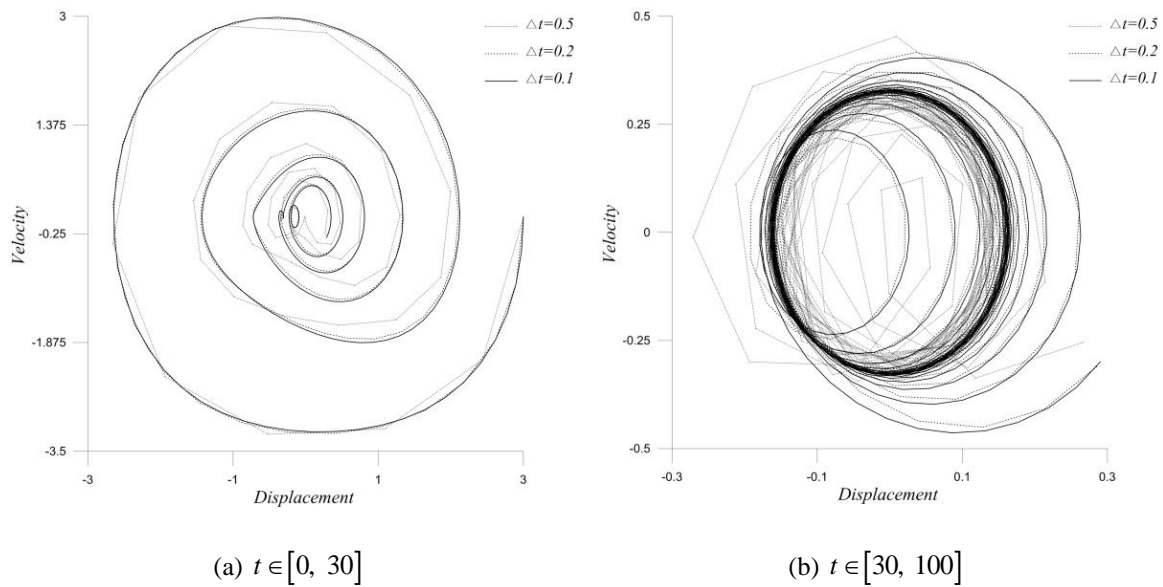

(a) $t \in [0, 30]$          (b) $t \in [30, 100]$

**Fig. 18. Phase portraits of T1E1 using three different time steps in the present method**

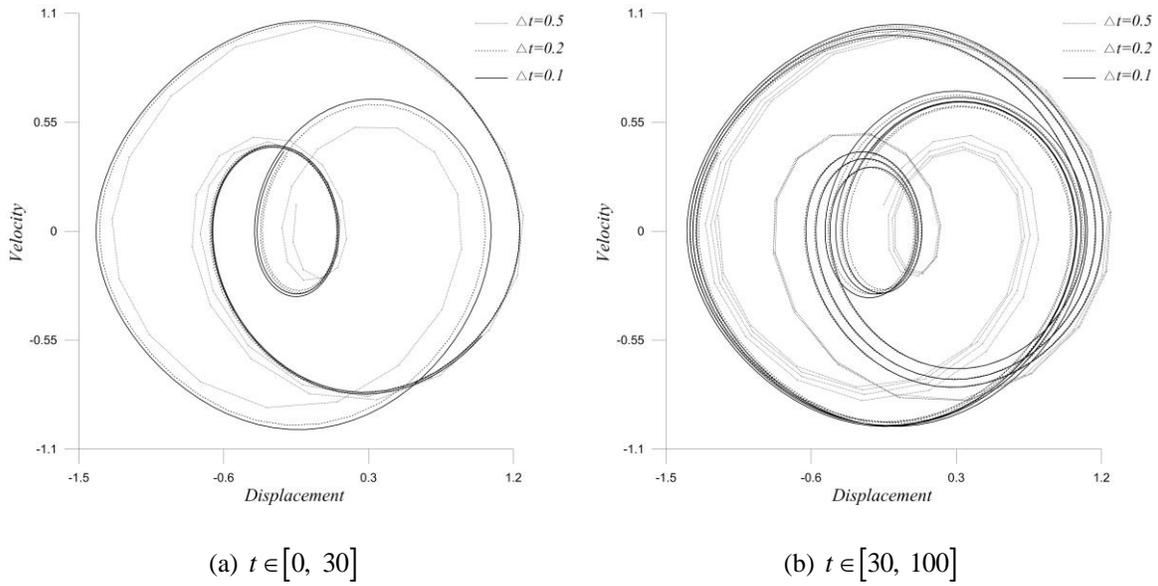

(a) $t \in [0, 30]$          (b) $t \in [30, 100]$

**Fig. 19. Phase portraits of T2E4 using three different time steps in the present method**

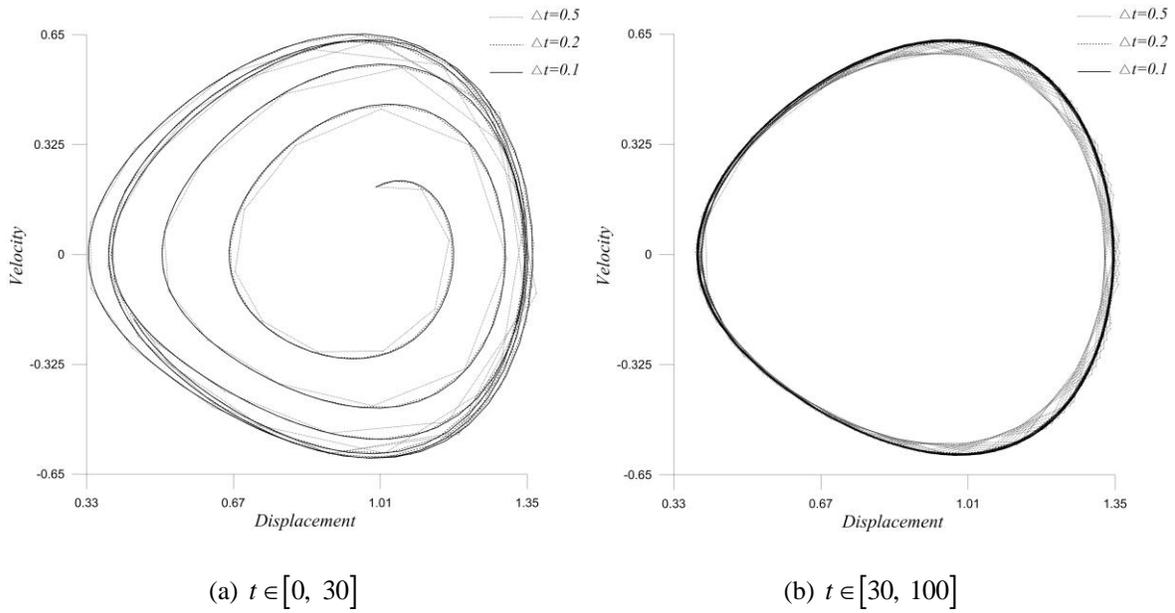

(a) $t \in [0, 30]$          (b) $t \in [30, 100]$

**Fig. 20. Phase portraits of T3E1 using three different time steps in the present method**

These figures illustrate convergence characteristics of the present method: as the time-step size decreases, simulation results get converged.

## 5. Conclusions

In this paper, based upon the extended framework of Hamilton's principle (EHP), a new variational formulation for the force-damped Duffing equation is presented. This formulation recovers all the governing differential equations as its Euler-Lagrange equation along with proper use of initial conditions. Thus, it provides elegant structure for the development of novel computational methods involving finite element presentation over time domain. To illustrate viability of temporal finite element approach stemming from this formulation, we also present the simplest temporal finite element method. The present numerical method adopts non-iterative algorithm and this method is validated with some examples in three types of the Duffing oscillation: hard spring; soft spring; inverted Duffing oscillation. For every example, we show that the present method with small time-step sizes yields almost identical results to RKF45 with default values of iteration and error tolerance in MAPLE software. Furthermore, with additional simulation for the selected representative example in each type of Duffing oscillation, we check that the present method has desirable convergence characteristics.

The present work only concerns the Duffing equation with cubic order nonlinearity, however, the simplistic beauty of the EHP allows us to further explore fractional or higher order nonlinearities. For example, the cubic-quintic nonlinearity ($\beta u^3 + \gamma u^5$) can be considered within the EHP by adding the term $\left(-\gamma / 6\right) u^6$ to Eq.(3). Also, there is still much room for the development of refined temporal finite element methods. In particular, future work will be directed toward development of higher-order time stepping method and varying time-step method.

# Acknowledgment

This research was supported by Basic Science Research Program through the National Research Foundation of Korea (NRF) funded by the Ministry of Science, ICT & Future Planning (No.2015R1A5A1037548).